\documentclass[12pt,twoside]{article} 
\usepackage[cp1251]{inputenc}
\usepackage[russian]{babel}
\usepackage{amsmath,amsfonts,amssymb}
\usepackage{latexsym,hhline,graphics,epsfig}

\newcounter{theorem}
\newcommand{\theor}{\par\refstepcounter{theorem}%
{\bf Теорема \thetheorem }.\,\,}

\newcounter{lemma}
\begin{document}

\pagestyle{myheadings} 
\markboth {\underline{\mbox{\!\!}\hspace{6,2cm} \small И. В.
Денега}}{\underline{\hspace{0.8cm}\small Обобщение некоторых
экстремальных задач о неналегающих областях... \hfill}}

\begin{flushleft}

\noindent {\small \textbf{УДК}\ \ 517.54}

\bigskip

{\bf И. В. Денега} (Институт математики НАН Украины, Киев)

\vskip 4mm

{\bf I. V. Denega} (Institute of mathematics of NAS of Ukraine,
Kyiv)

\vskip 4mm
{\bf Обобщение некоторых экстремальных задач о \\[1mm]неналегающих областях со свободными
 полюсами}

\vskip 4mm
{\bf Generalization of some extremal problems of \\[1mm]non-overlapping domains with free poles}
\end{flushleft}


Some results related to extremal problems with free poles are
generalized. The date have been obtained applying the known methods,
which are described in an earlier studies. Sufficiently good
numerical results of $\gamma$ were obtained.

В роботі узагальнені деякі екстремальні задачі про області, що не
перетинаються з вільними полюсами. Отримано покращене числове
значення $\gamma$.

В работе обобщены некоторые результаты экстремальных задач о
неналегающих областях со свободными полюсами. Получено улучшеное
числовое значение $\gamma$.

\thispagestyle{empty}
\newpage
\textbf{1. Введение.} В геометрической теории функций комплексного
переменного экстремальные задачи о неналегающих областях являются
хорошо известным классическим направлением. Под задачами такого рода
мы понимаем определение максимума произведения внутренних радиусов
попарно неналегающих областей, удовлетворяющих определенным
условиям. Возникновение данного направления геометрической теории
функций комплексной переменной связано с классической работой М.А.
Лаврентьева \cite{1}, в которой, в частности, была впервые
поставлена и решена задача о максимуме произведения конформных
радиусов двух непересекающихся односвязных областей. В дальнейшем
тематика связаная с изучением задач о неналегающих областях получила
развитие в работах  [1 -- 12]. В данной работе обобщены некоторые
результаты, полученные в \cite{8}.

Пусть $\mathbb{N}$, $\mathbb{R}$ -- множество натуральных и
вещественных чисел соответственно, $\mathbb{C}$ -- комплексная
плоскость, $\overline{\mathbb{C}}=\mathbb{C}\bigcup\{\infty\}$ --
одноточечная компактификация и $\mathbb{R^{+}}=[0,\infty)$.

Пусть $r(B,a)$ -- внутренний радиус области $B\subset\overline{C}$,
относительно точки $a\in B$ (см. напр. \cite{5}) и
$\chi(t)=\frac{1}{2}(t+t^{-1})$.

Пусть $n\in \mathbb{N}$. Систему точек $A_{n}:=\left\{a_{k} \in
\mathbb{C}:\, k=\overline{1,n}\right\},$ назовем \textbf{\emph{$n$ -
лучевой}}, если $|a_{k}|\in\mathbb{R^{+}}$ при $k=\overline{1,n}$,
$0=\arg a_{1}<\arg a_{2}<\ldots<\arg a_{n}<2\pi$. Обозначим при этом
$$\theta_{k}:=\arg a_{k},\,a_{n+1}:=a_{1},\, \theta_{n+1}:=2\pi,$$
$$\alpha_{k}:=\frac{1}{\pi}\arg \frac{a_{k+1}}{a_{k}},\, \alpha_{n+1}:=\alpha_{1},\,k=\overline{1, n}.$$

Для произвольной $n$-лучевой системы точек $A_{n}=\{a_{k}\}$ и
$\gamma\in\mathbb{R^{+}}$ полагаем
$$\mathcal{L}^{(\gamma)}(A_n):=\prod\limits_{k=1}^n\left[
\chi\left(\Bigl|\frac{a_k}{a_{k+1}}\Bigr|^\frac{1}{2\alpha_k}\right)\right]^{1-\frac{1}{2}\gamma\alpha_k^2}\cdot
\prod\limits_{k=1}^n|a_k|^{1+\frac{1}{4}\gamma(\alpha_k+\alpha_{k-1})}.$$

Целью данной работы является получение точных оценок сверху для
функционалов следующего вида

\begin{equation}\label{1a}J_{\gamma}=r^\gamma\left(B_0,0\right)\prod\limits_{k=1}^n r\left(B_k,a_k\right),\end{equation}
\begin{equation}\label{1b}
I_{\gamma}=\left[r\left(B_0,0\right)r\left(B_\infty,\infty\right)\right]^{\gamma}\prod\limits_{k=1}^n
r\left(B_k,a_k\right),\end{equation} где $\gamma\in\mathbb{R^{+}}$,
$A_{n}=\{a_{k}\}_{k=1}^{n}$ -- $n$-лучевая система точек, $a_{0}=0,$
$\{B_{k}\}_{k=0}^{n}$ -- система неналегающих областей (то есть
$B_{p}\cap B_{j}={\O}$ при $p\neq j$) таких, что $a_{k}\in B_{k}$
при $k=\overline{0, n}$.

\textbf{2. Основные результаты.} \textbf{\theor}{\it Пусть $n\in
\mathbb{N}$, $n\geqslant 2$ и $\gamma\in (0, 1]$. Тогда для любой
$n$-лучевой системы точек $A_n=\{a_k\}_{k=1}^n$,
$\mathcal{L}^{(\gamma)}\left(A_n\right)=1$ и любого набора взаимно
непересекающихся областей $B_k$, $a_k\in
B_k\subset\overline{\mathbb{C}}$, $k=\overline{0,n}$, справедливо
неравенство
\begin{equation}\label{2a}J_{\gamma}\leqslant\frac{4^{n+\frac{\gamma}{n}}
\gamma^\frac{\gamma}{n}n^n}{(n^2-\gamma)^{n+\frac{\gamma}{n}}}
\left(\frac{n-\sqrt{\gamma}}{n+\sqrt{\gamma}}\right)^{2\sqrt{\gamma}}.\end{equation}
Знак равенства в этом неравенстве достигается, когда $a_k$ и $B_k$,
$k=\overline{0,n}$, являются, соответственно, полюсами и круговыми
областями квадратичного дифференциала}
\begin{equation}\label{3a}Q(w)dw^2=-\frac{(n^2-\gamma)w^n+\gamma}{w^2(w^n-1)^2}\,dw^2.\end{equation}

\textbf{\theor}{\it Пусть $n\in \mathbb{N}$, $n\geqslant 2$ и
$\gamma=\frac{1}{2}$. Тогда для любой $n$-лучевой системы точек
$A_n=\{a_k\}_{k=1}^n$, $\mathcal{L}^{(0)}\left(A_n\right)=1$ и
любого набора взаимно непересекающихся областей $B_k$, $a_k\in
B_k\subset\overline{\mathbb{C}}$, $k=\overline{0,n}$, справедливо
неравенство
\begin{equation}\label{12a}
\left[r\left(B_0,0\right)r\left(B_\infty,\infty\right)\right]^{\frac{1}{2}}\prod\limits_{k=1}^n
r\left(B_k,a_k\right)\leqslant\frac{2^{2n+\frac{1}{n}}}{(n^2-2)^{\frac{1}{n}+\frac{n}{2}}}
\left(\frac{n-\sqrt{2}}{n+\sqrt{2}}\right)^{\sqrt{2}}.
\end{equation}
Знак равенства в этом неравенстве достигается, когда $a_k$ и $B_k$
являются, соответственно, полюсами и круговыми областями
квадратичного дифференциала}
\begin{equation}\label{13a}Q(w)dw^2=-\frac{w^{2n}+w^{n}(2n^2-2)+1}{w^2(w^n-1)^2}\,dw^2.\end{equation}

\textbf{\textit{Доказательство теоремы 1.}} Совершим разделяющее
преобразование системы областей $\{B_{k}\}_{k=1}^{n}$. Положим
$$E_k:=E_k(A_n):=\{w\in\mathbb{C}\backslash:\,\theta_k<\arg w<\theta_{k+1}\}.$$
Рассмотрим функцию
$\zeta=\pi_k(w)=-i\left(e^{-i\theta_k}w\right)^\frac{1}{\alpha_k}$,\quad
$k=\overline{1,n}$. При каждом $k=\overline{1,n}$ зафиксируем ту
ветвь многозначной аналитической функции $\pi_k(w)$, которая
осуществляет однолистное и конформное отображение $E_k$ на правую
полуплоскость $\text{Re}\,\zeta>0$.

Пусть $\Omega_k^{(1)}$, $k=\overline{1,n}$, обозначает область
плоскости $\zeta$, полученную в результате объединения связной
компоненты множества $\pi_k(B_k\bigcap\overline{E}_k)$, содержащей
точку $\pi_k(a_k)$, со своим симметричным отражением относительно
мнимой оси. В свою очередь, через $\Omega_k^{(2)}$,
$k=\overline{1,n}$, обозначаем область плоскости $\mathbb{C}_\zeta$,
полученную в результате объединения связной компоненты множества
$\pi_k(B_{k+1}\bigcap\overline{E}_k)$, содержащей точку
$\pi_k(a_k)$, со своим симметричным отражением относительно мнимой
оси, $B_{n+1}:=B_1$, $\pi_n(a_{n+1}):=\pi_n(a_1)$. Кроме того,
$\Omega_k^{(0)}$ будет обозначать область плоскости
$\mathbb{C}_\zeta$ , полученную в результате объединения связной
компоненты множества $\pi_k(B_0\bigcap\overline{E}_k)$, содержащей
точку $\zeta=0$, со своим симметричным отражением относительно
мнимой оси. Обозначим $\pi_k(a_k):=\omega_k^{(1)}$,
$\pi_k(a_{k+1}):=\omega_k^{(2)}$, $k=\overline{1,n}$,
$\pi_n(a_{n+1}):=\omega_n^{(2)}$.

Из определения функций $\pi_k$ вытекает, что
$$|\pi_k(w)-\omega_k^{(1)}|\sim\frac{1}{\alpha_k}|a_k|^{\frac{1}{\alpha_k}-1}\cdot|w-a_k|,\quad
w\rightarrow a_k,\quad w\in\overline{E_k},$$
$$|\pi_k(w)-\omega_k^{(2)}|\sim\frac{1}{\alpha_k}|a_{k+1}|^{\frac{1}{\alpha_k}-1}\cdot|w-a_{k+1}|,\quad
w\rightarrow a_{k+1},\quad w\in\overline{E_k},$$
$$|\pi_k(w)|\sim|w|^\frac{1}{\alpha_k},\quad
w\rightarrow 0,\quad w\in\overline{E_k},$$ Тогда, используя
соответствующие результаты работ [4--6], получаем неравенства
\begin{equation}\label{4a}r\left(B_k,a_k\right)\leqslant\left[\frac{r\left(\Omega_k^{(1)},\omega_k^{(1)}\right)
\cdot
r\left(\Omega_k^{(2)},\omega_k^{(2)}\right)}{\frac{1}{\alpha_k}
|a_k|^{\frac{1}{\alpha_k}-1}\cdot\frac{1}{\alpha_{k-1}}
|a_k|^{\frac{1}{\alpha_{k-1}}-1}}\right]^\frac{1}{2},\end{equation}
$$k=\overline{1,n},\quad \Omega_0^{(2)}:=\Omega_n^{(2)},\quad \omega_0^{(2)}:=\omega_n^{(2)},$$
\begin{equation}\label{5a}r\left(B_0,0\right)\leqslant\left[\prod \limits_{k=1}^n
r^{\alpha_k^2}\left(\Omega_k^{(0)},0\right)\right]^\frac{1}{2}.\end{equation}
Отсюда следует оценка для исследуемого функционала (\ref{1a})
$$J_{\gamma}\leqslant
\left[\prod\limits_{k=1}^n
r^{\gamma\alpha_k^2}\left(\Omega_k^{(0)},0\right)
\prod\limits_{k=1}^n r\left(\Omega_k^{(1)},\omega_k^{(1)}\right)
r\left(\Omega_k^{(2)},\omega_k^{(2)}\right)\right]^\frac{1}{2}\times$$
$$\times
\prod\limits_{k=1}^n\alpha_k\cdot \prod\limits_{k=1}^n
\frac{|a_k|}{\left[|a_k
a_{k+1}|^\frac{1}{\alpha_k}\right]^\frac{1}{2}}=\prod\limits_{k=1}^n\alpha_k\cdot
\prod\limits_{k=1}^n \frac{|a_k|}{|a_k
a_{k+1}|^\frac{1}{2\alpha_k}}\times$$
\begin{equation}\label{6a}\times \left[\prod\limits_{k=1}^n
r^{\gamma\alpha_k^2}\left(\Omega_k^{(0)},0\right)
\prod\limits_{k=1}^n r\left(\Omega_k^{(1)},\omega_k^{(1)}\right)
r\left(\Omega_k^{(2)},\omega_k^{(2)}\right)\right]^\frac{1}{2}.\end{equation}

Выражение (\ref{6a}), стоящее в скобках последней формулы,
представляет собой произведение значений функционала
$r^{\beta^2}(\Omega_k^{(0)},0)r(\Omega_k^{(1)},\omega_k^{(1)})r(\Omega_k^{(2)},\omega_k^{(2)})$
на тройках неналегающих областей
$\left(\Omega_k^{(0)},\Omega_k^{(1)},\Omega_k^{(2)}\right)$
плоскости $\zeta$.

Перейдем далее к инвариантной форме этого функционала, что дает
существенные преимущества. С учетом результатов работы \cite{7},
получаем
$$r^\gamma\left(B_0,0\right)\prod\limits_{k=1}^n r\left(B_k,
a_k\right)\leqslant
\left(\prod\limits_{k=1}^n\alpha_k\right)\cdot\prod\limits_{k=1}^n
\frac{|a_k|}{|a_k a_{k+1}|^\frac{1}{2\alpha_k}}\times$$
$$\times \left\{\prod\limits_{k=1}^n\frac{
r^{\gamma\alpha_k^2}\left(\Omega_k^{(0)},0\right)\cdot
r\left(\Omega_k^{(1)},\omega_k^{(1)}\right)\cdot
r\left(\Omega_k^{(2)},\omega_k^{(2)}\right)}{|\omega_k^{(1)}\cdot
\omega_k^{(2)}|^{\gamma\alpha_k^2}|\omega_k^{(1)}-
\omega_k^{(2)}|^{2-\gamma\alpha_k^2}}\right\}^\frac{1}{2}\times$$
$$\times\left[\prod\limits_{k=1}^n |\omega_k^{(1)}\cdot
\omega_k^{(2)}|^{\gamma\alpha_k^2}|\omega_k^{(1)}-
\omega_k^{(2)}|^{2-\gamma\alpha_k^2}\right]^\frac{1}{2}=$$
$$=\left(\prod\limits_{k=1}^n\alpha_k\right)\cdot\prod\limits_{k=1}^n
\frac{|a_k|}{|a_k a_{k+1}|^\frac{1}{2\alpha_k}}\times$$
$$\times\left(\prod\limits_{k=1}^n|\omega_k^{(1)}-
\omega_k^{(2)}|\right)\left(\prod\limits_{k=1}^n\frac{|\omega_k^{(1)}\cdot
\omega_k^{(2)|}}{|\omega_k^{(1)}-
\omega_k^{(2)}|}\right)^{\frac{\gamma\alpha_{k}^{2}}{2}}\times$$
$$\times \left\{\prod\limits_{k=1}^n\frac{
r^{\gamma\alpha_k^2}\left(\Omega_k^{(0)},0\right)\cdot
r\left(\Omega_k^{(1)},\omega_k^{(1)}\right)\cdot
r\left(\Omega_k^{(2)},\omega_k^{(2)}\right)}{|\omega_k^{(1)}\cdot
\omega_k^{(2)}|^{\gamma\alpha_k^2}|\omega_k^{(1)}-
\omega_k^{(2)}|^{2-\gamma\alpha_k^2}}\right\}^\frac{1}{2}=$$
\begin{equation}\label{7a}\end{equation}
$$=2^{n}\cdot\left(\prod\limits_{k=1}^n\alpha_k\right)\cdot\prod\limits_{k=1}^n
\chi\left(\left|\frac{a_k}{a_{k+1}}\right|^\frac{1}{2\alpha_k}\right)|a_k|\times$$
$$\times2^{-\frac{\alpha}{2}\sum\limits_{k=1}^n\alpha_k}\left[\prod\limits_{k=1}^n
\chi\left(\left|\frac{a_k}{a_{k+1}}\right|^\frac{1}{2\alpha_k}\right)\right]^{-\frac{\gamma\alpha_{k}^{2}}{2}}
\left(\prod\limits_{k=1}^n\left|\frac{a_{k+1}}{a_{k}}\right|\right)^{\frac{\gamma\alpha_{k}^{2}}{2}}\times$$
$$\times \left\{\prod\limits_{k=1}^n\frac{
r^{\gamma\alpha_k^2}\left(\Omega_k^{(0)},0\right)\cdot
r\left(\Omega_k^{(1)},\omega_k^{(1)}\right)\cdot
r\left(\Omega_k^{(2)},\omega_k^{(2)}\right)}{|\omega_k^{(1)}\cdot
\omega_k^{(2)}|^{\gamma\alpha_k^2}|\omega_k^{(1)}-
\omega_k^{(2)}|^{2-\gamma\alpha_k^2}}\right\}^\frac{1}{2}=$$
$$=2^{n-\frac{\gamma}{2}\sum\limits_{k=1}^n\alpha_k^2}\cdot
\left(\prod\limits_{k=1}^n\alpha_k\right)
\cdot\prod\limits_{k=1}^n\left[\chi\left(\Bigl|\frac{a_k}{a_{k+1}}\Bigr|
^\frac{1}{2\alpha_k}\right)\right]^{1-\frac{\gamma\alpha_k^2}{2}}\times$$
$$\times\prod\limits_{k=1}^n|a_k|^{1+\frac{1}{4}\gamma(\alpha_k+\alpha_{k-1})}\times$$
$$\times\left\{\prod\limits_{k=1}^n\frac{
r^{\gamma\alpha_k^2}\left(\Omega_k^{(0)},0\right)\cdot
r\left(\Omega_k^{(1)},\omega_k^{(1)}\right)\cdot
r\left(\Omega_k^{(2)},\omega_k^{(2)}\right)}{|\omega_k^{(1)}\cdot
\omega_k^{(2)}|^{\gamma\alpha_k^2}|\omega_k^{(1)}-
\omega_k^{(2)}|^{2-\gamma\alpha_k^2}}\right\}^\frac{1}{2}.$$ Здесь
использованы соотношения
$|\omega_k^{(1)}|=|a_k|^\frac{1}{\alpha_k}$,
$|\omega_k^{(2)}|=|a_{k+1}|^\frac{1}{\alpha_k}$,
$|\omega_k^{(1)}-\omega_k^{(2)}|=|a_k|^\frac{1}{\alpha_k}+|a_{k+1}|^\frac{1}{\alpha_k}$.
Функционал, стоящий в фигурных скобках в выражении (\ref{7a}),
инвариантен при любых конформных автоморфизмах
$\overline{\mathbb{C}}$.

При каждом $k=\overline{1,n}$ несложно указать конформный
автоморфизм $\zeta=T_k(z)$ плоскости комплексных чисел
$\overline{\mathbb{C}}$ такой, что $T_k(0)=0$,\,
$T_k\left(g_k^{(s)}\right)=(-1)^s\cdot i$,\,
$\Omega_k^{(q)}:=T_k\left(G_k^{(q)}\right)$,\, $k=\overline{1,n}$,\,
$s=1,2$,\, $q=0,1,2$. Инвариантность относительно конформных
автоморфизмов $\overline{\mathbb{C}}$ функционала
$$J_3(\alpha_1,\alpha_2,\alpha_3)=
\frac{r^{\alpha_1}(B_1,a_1)\cdot r^{\alpha_2}(B_2,a_2)\cdot
r^{\alpha_3}(B_3,a_3)}
{|a_1-a_2|^{\alpha_1+\alpha_2-\alpha_3}\cdot|a_1-a_3|^{\alpha_1-\alpha_2+\alpha_3}\cdot
|a_2-a_3|^{-\alpha_1+\alpha_2+\alpha_3}},$$
$\alpha_k\in\mathbb{R}^+$,\, $a_k\in
B_k\subset\overline{\mathbb{C}}$,\, $B_k\bigcap B_p=\varnothing$,\,
$k=1,2,3$,\, $p=1,2,3$,\, $k\neq p$, по-видимому, впервые указана в
работе \cite{7}. С учетом этой инвариантности, получаем
$$J_{\gamma}\leqslant
2^{n-\frac{\gamma}{2}\sum\limits_{k=1}^n\alpha_k^2}\cdot
\left(\prod\limits_{k=1}^n\alpha_k\right)
\cdot\prod\limits_{k=1}^n\left[\chi\left(\Bigl|\frac{a_k}{a_{k+1}}\Bigr|
^\frac{1}{2\alpha_k}\right)\right]^{1-\frac{\gamma\alpha_k^2}{2}}
\prod\limits_{k=1}^n|a_k|^{1+\frac{1}{4}\gamma(\alpha_k+\alpha_{k-1})}\times$$

$$\times\prod\limits_{k=1}^n\left\{\frac{r^{\alpha_k^2\gamma}\left(\Omega_k^{(0)},0\right)\cdot
r\left(\Omega_k^{(1)},-i\right)\cdot
r\left(\Omega_k^{(2)},i\right)}{2^{2-\gamma\alpha_{k}^{2}}}\right\}^\frac{1}{2}=$$
$$=2^{n-\frac{\gamma}{2}\sum\limits_{k=1}^n\alpha_k^2}\left(\prod\limits_{k=1}^n\alpha_k\right)\cdot\mathcal{L}^{(\gamma)}(A_n)\cdot2^{-n+\frac{\gamma}{2}\sum\limits_{k=1}^n\alpha_k^2}\times$$
$$\times\left[\prod\limits_{k=1}^n
r^{\alpha_k^2\gamma}\left(\Omega_k^{(0)},0\right)\cdot
r\left(\Omega_k^{(1)},-i\right)\cdot
r\left(\Omega_k^{(2)},i\right)\right]^\frac{1}{2}\leqslant$$
$$\leqslant\left(\prod\limits_{k=1}^n\alpha_k\right)\cdot\mathcal{L}^{(\gamma)}(A_n)\cdot
\left[\prod\limits_{k=1}^n
r^{\alpha_k^2\gamma}\left(\Omega_k^{(0)},0\right)\cdot
r\left(\Omega_k^{(1)},-i\right)\cdot
r\left(\Omega_k^{(2)},i\right)\right]^\frac{1}{2}.$$

В результате проведенных вычислений исходная задача сведена к оценке
сверху функционала $r^{\beta^2}(B_0,0)r(B_1,i)r(B_2,-i)$ на классе
троек попарно непересекающихся областей $\{B_0,B_1,B_{2}\}$ таких,
что $0\in B_0$, $i\in B_1$, $-i\in B_2$,
$B_k\subset\overline{\mathbb{C}}$, $k=0,1,2$.

Далее используем метод, предложенный В.Н. Дубининым при
доказательстве теоремы 4 \cite{4}, и учитывая результат Л.И.
Колбиной \cite{7}, получаем
\begin{equation}\label{9a}
r^{\beta^2}(B_0,0)r(B_1,i)r(B_2,-i)\leqslant\Psi(\beta)=\end{equation}
$$=2^{\beta^2+6}\cdot\beta^{\beta^2+2}(2-\beta)^{-\frac{1}{2}(2-\beta)^2}\cdot
(2+\beta)^{-\frac{1}{2}(2+\beta)^2},\quad\beta\in[0,2].$$

Рассмотрим экстремальную задачу:
\begin{equation}\label{10a}\prod\limits_{k=1}^n\Psi(\beta)\rightarrow\sup;\quad \sum\limits_{k=1}^n\beta_{k}=2.
\end{equation}
Необходимые условия имеют вид
\begin{equation}\label{11a}\frac{\Psi'(\beta)}{\Psi(\beta)}=\frac{-\lambda}{\prod\limits_{k=1}^n\Psi(\beta)},\quad
k=\overline{1,n}.\end{equation}

Покажем, что все $\beta_{k}$ равны между собой. Исследуем поведение
функции
$F(\beta)=\frac{\Psi'(\beta)}{\Psi(\beta)}=2\beta\ln(2\beta)+\frac{2}{\beta}+(2-\beta)\ln(2-\beta)-(2+\beta)\ln(2+\beta)$
на промежутке $\beta\in[0,2]$. Она строго убывает на промежутке $(0;
\beta_{0}]$, $\beta_{0}\in (1,32; 1,33)$ и возрастает на
$[\beta_{0}; 2)$. Введем вспомогательную функцию
$\mathfrak{F}(\beta)=F(\beta)-F(2-\beta), \:0<\beta\leq1.$
$\mathfrak{F}(\beta)$ положительна на интервале $(0, 1)$, поскольку
$F(2-\beta)$ лежит ниже чем $F(\beta)$. Если теперь хотя бы одно из
$\beta_{k}$ больше 1, например, $\beta_{k'}>1$, то для остальных
$\beta_{k}$ имеем неравенства
$\beta_{k}\leq2-\beta_{k'}<1<\beta_{0}$ и потому $F(\beta_{k})\geq
F(2-\beta_{k'})>F(2-(2-\beta_{k'}))=F(\beta_{k'}).$ Получили
противоречие с (\ref{11a}). Поэтому для произвольного $k$
справедливо условие $\beta_{k}\in (0, 1]$ и в силу равенств
(\ref{10a}), а также монотонности $F(\beta)$ на $(0, \beta_{0}]$ в
точке предполагаемого экстремума, имеем $\beta_{k}=\frac{2}{n},$
$k=\overline{1,n}.$ Отсюда следует, что точка $(\frac{2}{n}, \ldots,
\frac{2}{n})$ действительно есть решением задачи (\ref{10a}).
Используя оценки (\ref{4a}), (\ref{5a}) и (\ref{9a}) получаем
доказательство теоремы.

Из свойств разделяющего преобразования \cite{4} получаем, что знак
равенства в неравенстве (\ref{2a}) достигается, когда точки $a_k$ и
области $B_k$, $k=\overline{0,n}$, являются, соответственно,
полюсами и круговыми областями квадратичного дифференциала
(\ref{3a}). Теорема 1 доказана.

В. Н. Дубинин доказал этот результат при $\gamma=1$ для любых
различных точек $a_k$, лежащих на окружности $|z|=1$ и любых попарно
непересекающихся областей $B_k$ (см. \cite{4,10}).

\textbf{Доказательство теоремы 2.} Мы сохраняем все обозначения для
разделяющего преобразования областей, введенные при доказательстве
теоремы 1 для областей $B_k$, $k=\overline{0,n}$. Кроме того,
$\Omega_k^{(\infty)}$ будет обозначать область плоскости
$\mathbb{C}_\zeta$ , полученную в результате объединения связной
компоненты множества $\pi_k(B_\infty\bigcap\overline{E}_k)$,
содержащей точку $\zeta=\infty$, со своим симметричным отражением
относительно мнимой оси. Семейство
$\left\{\Omega_k^{(\infty)}\right\}_{k=1}^n$, является результатом
разделяющего преобразования произвольной области $B_\infty$,
$\infty\in B_\infty\subset\overline{\mathbb{C}}$, относительно
семейств $\left\{E_k\right\}_{k=1}^n$ и
$\left\{\pi_k\right\}_{k=1}^n$ в точке $\zeta=\infty$.

По теореме 2 \cite{4} имеем
\begin{equation}\label{15a}
r(B_\infty,\infty)\leq\left[\prod \limits_{k=1}^n
r^{\alpha_k^2}\left(\Omega_k^{(\infty)},\infty\right)\right]^\frac{1}{2}.
\end{equation}

Используя оценки (\ref{4a}), (\ref{5a}), (\ref{15a}), получаем
$$\left[r\left(B_0,0\right)r\left(B_\infty,\infty\right)\right]^{\frac{1}{2}}\prod\limits_{k=1}^n
r\left(B_k,a_k\right)\leqslant2^{n}\cdot\left(\prod\limits_{k=1}^n\alpha_k\right)\times$$
$$\times\mathcal{L}^{(0)}(A_n)\cdot
\left[\prod\limits_{k=1}^n
r^{\alpha_k^2\gamma}\left(\Omega_k^{(0)},0\right)\cdot
r^{\alpha_k^2\gamma}\left(\Omega_k^{(\infty)},\infty\right)\cdot
r\left(\Omega_k^{(1)},-i\right)\cdot
r\left(\Omega_k^{(2)},i\right)\right]^\frac{1}{2}.$$

Теорема 6 работы \cite{4} дает
\begin{equation}\label{14a}
\left[r(B_0,0)r(B_\infty,\infty)\right]^{\frac{1}{2}}r(B_1,i)r(B_2,-i)\leqslant\Psi(\beta)=\end{equation}
$$=8\beta^{2\beta^2+2}\cdot|1-\beta|^{-(1-\beta)^2}\cdot(1+\beta)^{-(1+\beta)^2},
\quad 0<\beta\leq\sqrt{2}.$$ Неравенство (\ref{14a}) получено В.Н.
Дубининым с использованием результатов Л.И. Колбиной \cite{7}.

Рассмотрим экстремальную задачу:
$$\prod\limits_{k=1}^n\Psi(\beta_{k})\rightarrow\sup;\quad
\sum\limits_{k=1}^n\beta_{k}=\sqrt{2}. $$

Введем функцию $F(\beta)=\frac{\Psi'(\beta)}{\Psi(\beta)}$.
Вычисления показывают, что эта функция убывает на промежутке $(0;
\beta_{0}]$ и возрастает на $[\beta_{0}; \sqrt{2})$\quad
$0,85<\beta_{0}<1$. Учитывая также, что $F(0,564)>0$,
$F(\sqrt{2})<0$, получаем, что разность $F(\beta)-F(\sqrt{2}-\beta)$
положительна на промежутке $0<\beta_{0}<\frac{\sqrt{2}}{2}$.
Аналогично доказательству теоремы 1, убеждаемся, что единственным
решением экстремальной задачи является точка $(\frac{\sqrt{2}}{n},
\ldots, \frac{\sqrt{2}}{n})$. Оценки (\ref{4a}), (\ref{5a}),
(\ref{15a}), (\ref{14a}) дают неравенство теоремы 2. Случай
равенства проверяется непосредственно. Теорема 2 доказана.

\emph{В заключение выражаю благодарность А. К. Бахтину за
постановку задач и полезные обсуждения.}


{\small
\begin{thebibliography}{99}
\bibitem{1} Лаврентьев М. А. {\sl К теории
конформных отображений}// Тр. Физ.-мат. ин-та АН СССР. -- 1934.--
5.-- С. 159 -- 245.
\bibitem{2} Голузин Г. М. {\sl Геометрическая теория функций комплексного
переменного}. -- М: Наука, 1966. -- 628 с.
\bibitem{heim} Хейман В К. {\sl Многолистные функции. - М.: Изд-во иностр.
лит}., 1960. -- 180 с.
\bibitem{3} Дубинин В. Н. {\sl Метод симметризации в задачах о неналегающих
областях} // Мат. сб. -- 1985. -- {\bf 128}, №~1. -- С.~110~--~123.
\bibitem{4}  Дубинин В. Н. {\sl Разделяющее преобразование областей и
задачи об экстремальном разбиении} // Зап. науч. сем. Ленингр.
отд-ния Мат. ин-та АН СССР. -- 1988. -- 168. -- С. 48 -- 66.
\bibitem{5}  Дубинин В. Н. {\sl Метод симметризации в геометрической
теории функций комплексного переменного} // Успехи мат. наук. --
1994. -- 49, № 1(295). -- С. 3 -- 76.
\bibitem{6}  Дубинин В. Н. {\sl Асимптотика модуля вырождающегося
конденсатора и некоторые ее применения} // Зап. науч. сем. ПОМИ. --
1997. -- 237. -- С. 56 -- 73.
\bibitem{7} Колбина Л. И. {\sl Конформное отображение единичного круга на
неналегающие области} // Вестник Ленинград. ун-та. -- 1955. -- {\bf
5}. -- С.~37~--~43.
\bibitem{kovalev} Ковалев Л.В. {\sl К задаче об экстремальном разбиении со свободными полюсами на окружности.}
// Дальневосточный матем. сборник. -- 1996. -- {\bf 2}. -- С.~96~--~98.
\bibitem{8}  Бахтина Г. П., Бахтин А. К. {\sl Разделяющее преобразование и задачи о неналегающих областях}
// Збірник праці Ін-ту мат-ки НАН Укр. -- 2006. -- Т. 3., № 4, -- 273 -- 281 с.
\bibitem{9}  Бахтин А. К., Бахтина Г. П., Зелинский Ю. Б.
{\sl Тополого-алгебраические структуры  и геометрические методы в
комплексном анализе.} // Праці ін-ту мат-ки НАН Укр. -- 2008. -- 308
с.
\bibitem{10}  Дубинин В. Н.{\sl Емкости конденсаторов и симметризация в
геометрической теории функций комплексного переменного.}
// Владивосток "Дальнаука" ДВО РАН -- 2009. -- 390с.

\end{thebibliography}}
\end{document}